%% file: phase2.tex
\documentclass{article}
\usepackage{latexsym}
\usepackage{amssymb}
\usepackage{amsmath}
\usepackage{enumerate}
\usepackage[dvips]{graphicx}
\usepackage[latin1]{inputenc}

\newcommand\beginpf{\noindent \emph{Proof: }}
\newcommand\eprf{\hfill$\Box$}

\newcommand\arr{\mathbb{R}}

\newcommand\scrA{{\cal A}}
\newcommand\scrB{{\cal B}}
\newcommand\scrF{{\cal F}}
\newcommand\scrM{{\cal M}}
\newcommand\scrP{{\cal P}}
\newcommand\scrR{{\cal R}}
\newcommand\scrO{{\cal O}}

\newcommand\scrQ{{\cal Q}}

\newcommand\ol{\overline}
\newcommand\ul{\underline}

\newcommand\HD{\mathrm{HD}}
\newcommand\supp{\mathrm{supp}}

\newcommand\remark{\noindent \emph{Remark: }}

\newtheorem{thm}{Theorem}
\newtheorem{dfn}[thm]{Definition}
\newtheorem{lem}[thm]{Lemma}
\newtheorem{prop}[thm]{Proposition}

\begin{document}

\title{Renormalisation-induced phase transitions for unimodal maps }
\author{Neil Dobbs}
\thanks{The author was supported by the EU training network ``Conformal Structures and Dynamics''.}

\date{\today}
\maketitle
\begin{abstract}
The thermodynamical formalism is studied for renormalisable maps of the interval and the natural potential $-t \log|Df|$. Multiple and indeed infinitely many phase transitions at positive $t$ can occur for some quadratic maps. All unimodal quadratic maps with positive topological entropy exhibit a phase transition in the negative spectrum. 
\end{abstract}
\input{b2}

\bibliography{references}
\bibliographystyle{plain}

\end{document}

%% file: b2.tex
\section*{Introduction}

This brief note is dedicated to the presentation of a hitherto unstudied phenomenon: phase transitions between positive entropy regimes for the natural potential $-t \log |Df|$ in unimodal dynamics. 

In Proposition \ref{prop:infrenorm}, perhaps the most striking result, we show the existence of infinitely-renormalisable quadratic maps which admit an infinity of phase transitions. In Proposition \ref{prop:negt}, we show that all interesting quadratic maps exhibit a phase transition in the negative spectrum. In Proposition \ref{prop:CEpht}, we show that for each $n \geq 0$, there is a positive measure set of Collet-Eckmann parameters for which the maps admit exactly $n$ phase transitions in a neighbourhood of $[0,1]$. 

In related prior examples of phase transitions (e.g.\ \cite{BruinKeller:Fib}, \cite{BT:EquilibriumInterval}, \cite{MakStas:NegSpec}, \cite{MakStas:Nonrec}, \cite{PrellbergSlawny}), all equilibrium states to one side of the phase transition have had zero entropy. 

Sarig has interesting results (\cite{Sarig:CMP06},  \cite{Sarig:CMP01}, \cite{Sarig:CRAS00}) on phase transitions in the cadre of countable state Markov shifts; however these are not smooth dynamical systems. 

We consider real maps $f : I \to I$ of a closed interval $I$ and of class $C^2$ which are unimodal: $f$ has exactly two branches of monotonicity separated by the unique critical point $c$, and $f(\partial I) \to \partial I$. 
The basic family of unimodal maps to consider is the quadratic family: to each positive real parameter $a$ is assigned a quadratic map $f_a : x \to ax(1-x)$. For such maps, any $f_a$-invariant probability measure necessarily lies in the interval $[0,1]$, and for parameters $a \in (0,4]$, $f_a :[0,1] \to [0,1]$. 

We are interested in the natural potential $-t \log |Df|$, where $t \in \arr$ is the potential parameter and $Df$ denotes the derivative. The entropy of a measure we denote $h_\mu$ and the Lyapunov exponent $\chi_\mu = \int \log |Df| d\mu$. Let $\scrM(f)$ denote the class of ergodic, $f$-invariant probability measures.
 For $\mu \in \scrM(f)$, we set 
$$
\scrF(\mu, t) := h_\mu - t\chi_\mu.
$$
This quantity is related to minus the \emph{free energy}.
The \emph{(variational) pressure} (at parameter $t$) is then defined as 
$$
\scrP(t) := \scrP_f(t) := \sup_{\mu \in \scrM(f)} \scrF(\mu,t) =  \sup_{\mu \in \scrM(f)} (h_\mu - t \chi_\mu).
$$
See \cite{PrzRLSm} for a discussion of various definitions of pressure.
Provided that the critical orbits are non-periodic, the pressure will be finite.
We usually drop from the notation the dependence of the pressure function on the dynamical system $f : I \to I$ considered. 
\emph{Equilibrium measures}, or \emph{states}, are measures which realise the supremum. It will also be useful to define, for compact, forward-invariant subsets $X$ of $I$,
$$\scrM(f,X) := \{ \mu \in \scrM(f) : \supp(\mu) \subset X\}$$ 
and then 
$$
\scrP(t,X) := \sup_{\mu \in \scrM(f,X)} \scrF(\mu,t).
$$ 

Note that for fixed $\mu$, $\scrF(\mu,t)$ is linear in $t$. 
Convexity of the pressure function follows easily. 
Moreover, if all measures have non-negative Lyapunov exponents, the pressure function is decreasing. Any measure in $\scrM(f)$ with negative Lyapunov exponent is supported on an attracting periodic orbit if the critical point of $f$ is non-flat \cite{Przytycki:LyapNonneg}. In particular, in the absence of periodic attractors,  Lyapunov exponents of measures are non-negative. A critical point is said to be \emph{non-flat} if there exists a $C^2$ diffeomorphism $\phi: \arr \to \arr$ with $\phi(0) = c$ such that $f \circ \phi$ is a polynomial near the origin.

One expects the pressure function to be smooth, or actually analytic, almost everywhere. We say there is a \emph{phase transition} at $t$ if $\scrP(t)$ is \emph{not differentiable at $t$}. This is not the only possible definition. It is weaker than \emph{non-analyticity of $\scrP(t)$ at $t$} and differs from the existence of  \emph{more than one equilibrium state}. At the phase transitions we exhibit there will be two distinct equilibrium states.  	 
 Typically, at phase transitions, the equilibrium states will not vary continuously with the potential parameter. 

We shall study how renormalisations in unimodal dynamics can lead to phase transitions as the equilibrium states jump between different transitive sets. We exhibit examples of quadratic maps with an infinity of phase transitions at potentials between zero and one, and so where the equilibrium states have positive entropy.

Interesting recent results for thermodynamics of \emph{non-renormalisable maps} of the interval with respect to the natural potential are to be found in the works of Bruin and Todd \cite{BT:EquilibriumInterval} and of Pesin and Senti \cite{PesinSenti:Moscow}. 

We shall now define renormalisations and explain in a general way how they can give rise to phase transitions before providing details a little later.

A \emph{restrictive interval} $J$ of period $k>1$ is a closed non-degenerate interval, strictly contained in $I$ and containing the unique critical point $c$ with the following properties.  For all $i, j$ with $0 \leq i < j < k$, $f^i(J)$ and $f^j(J)$ overlap in at most one point; $f^k(\partial J) \in \partial J$ and $f^k(J) \subset J$. 

The map $f$ is renormalisable of type $k$ if $f$ has a (maximal) restrictive interval $J$ of minimal period $k$. The unimodal map $\scrR f := f^k_{|J} : J \to J$ is called a \emph{renormalisation} of $f$ of type $k$. If $k = 2$ the renormalisation is also called a Feigenbaum renormalisation. Since the renormalisation of $f$ is unimodal, its (topological) entropy is bounded by $\log 2$. Thus the entropy of $f$ restricted to $\bigcup_i f^i(J)$  is bounded by $(1/k)\log 2$.

The well known Feigenbaum map is the quadratic map which is infinitely renormalisable of type 2, that is, each renormalisation is of type 2. It has zero topological entropy. 

A unimodal map with positive topological entropy $\log s$ is renormalisable of type 2 \emph{if and only if} $\log s \leq (1/2)\log 2$. This can be seen by considering the semi-conjugacy with the piecewise linear \emph{tent map} with slope $\pm s$ --- see for example section III.4 of \cite{DeMeloVanStrien}, which indeed provides a good background for this paper.

If $f$ is (once-) renormalisable of type $k>2$ with restrictive interval $J$ of period $k$, then $f$ has topological entropy strictly greater than $(1/2)\log 2$, as noted before. Denote by $X_0$ the Cantor set of points which never enter the interior of $J$ under iteration by $f$. Since the entropy of $f$ restricted to the forward orbit of $J$ is less than or equal to $(1/2)\log 2$, the entropy of $f$ restricted to $X_0$ is necessarily the entropy of $f$.      Thus all measures of sufficiently large entropy are supported on $X_0$. 

Since, for $t$ close to zero, the entropy of any equilibrium state must have entropy close to the topological entropy of $f$, it follows that any such equilibrium state must be supported on $X_0$. 
Supposing that $X_0$ is hyperbolic,
existence of such equilibrium states follows from the general theory for hyperbolic maps (\cite{Ruelle:Book}). 
On the other hand, for sufficiently large $t>0$, one can expect any equilibrium measure to give positive measure to (the interior of) $J$.

We now explain why $X_0$ often is hyperbolic. 
The Mañé hyperbolicity theorem \cite{Mane:Hyperbolicity}, states that if $f$ is of class $C^2$, any forward invariant compact set disjoint from critical points and non-repelling periodic points  is a hyperbolic set.  
If $f$ has negative Schwarzian derivative, then any parabolic periodic point is attracting on one side and contains a (the) critical point in its immediate basin of attraction. In particular, if $\partial J$ does not contain a parabolic periodic point, then $X_0$ contains no parabolic periodic points, so $X_0$ is hyperbolic.

We would like to insist that the dynamics of $f$ restricted to $X_0$ and of $\scrR f : J \to J$ are essentially unconnected. If one no longer requires $f$ to be analytic then one can perturb the maps essentially independently. For this reason one can expect that (even for quadratic maps) a renormalisation of type $>2$ will lead to a phase transition. This is perhaps the key idea of the paper.

It follows that with strings of renormalisations one can find unimodal maps whose pressure functions admit various different behaviours: linear followed by strictly convex; strictly convex followed by linear followed by strictly convex followed by linear; piecewise linear etc. Constructing such maps is left as an exercise for the curious reader. 

One could ask whether similar phenomena to those exhibited in this paper can occur for rational maps. They cannot: in \cite{Me:Rat} we show that there exists at most one equilibrium state with positive entropy for each $t$. This is essentially due to the \emph{eventually onto} property in rational dynamics. 


\section*{Results}

We now proceed to formally state some results and outline their proofs, starting with some  simple observations. 
\begin{dfn}
Let $l_1 \leq l_2 \in \arr \cup \{\pm \infty\}$. We say a function $h:X \to \arr$, where $X \subset \arr$, has \emph{$(l_1, l_2)$-bounded slope} if, for all $x<y \in X$,
$$
l_1(y-x) \leq h(y) - h(x) \leq l_2(y-x).
$$
\end{dfn}
\begin{lem}\label{lem:bddslope}
Let $h, h' : \arr \to \arr$ be two real functions. Suppose $h$ has $(l_1, l_2)$-bounded slope and $h'$ has $(m_1, m_2)$-bounded slope, where $l_1 \leq l_2 < m_1 \leq m_2$. Then there exists a unique $x_0 \in \arr$ for which $h(x_0) = h'(x_0)$. Moreover, the function $\max(h,h')$ is not differentiable at $x_0$, and it coincides with $h$ on $\{x \leq x_0\}$ and with $h'$ on $\{x \geq x_0\}$.
\end{lem}
\beginpf
Evident.
\eprf

For $f :I \to I$ be a map of the interval and $X$ a compact, forward-invariant subset of $I$, let $\overline{\chi}(X)$ denote the supremum of the Lyapunov exponents of measures in $\scrM(f,X)$ and $\underline{\chi}(X)$ the infimum. 
\begin{lem}\label{lem:pslope}
The restricted pressure function $P(t,X)$ has $(-\overline{\chi}(X), -\underline{\chi}(X))$-bounded slope.
\end{lem}
\beginpf 
This follows easily from the definitions.
\eprf

Define $\scrQ$ as the set of $a \in [3,4]$ for which the critical point is not periodic. Note that all \emph{interesting} quadratic maps have parameters in $[3,4]$: if $a$ is greater than $4$ then $f$ restricted to the non-wandering set is hyperbolic and conjugate to the full shift on two symbols; for $a \in [0,3]$ there are only a finite number of periodic orbits. 

For $a\in [3,4] \setminus \scrQ$, the pressure, as defined, is infinite for all strictly positive $t$. One can alternatively restrict the definition of pressure to the supremum over measures living on the Julia set, and thus avoid this problem, since the Lyapunov exponent of a  measure being negative implies the measure sits on a periodic attractor (\cite{Przytycki:LyapNonneg}). 
The following proposition could then go, \emph{for all quadratic unimodal maps $f$ with positive topological entropy, there is a phase transition\ldots} In any case, for $a \notin \scrQ$, the reader can formulate corresponding statements without too much difficulty.
\begin{prop} \label{prop:negt} For all $a \in \scrQ$, the pressure function of the quadratic map $f_a$, with topological entropy denoted $h_a$, admits a phase transition at some $t \leq -h_a /(\log a - h_a)$.
\end{prop}
\beginpf 
The phase transition will be caused by the unit mass $\mu_0$ sitting on the repelling fixed point at zero, whose Lyapunov exponent is $\log a$. Thus $\scrP(t,\{0\}) = -t \log a$ and in particular $\scrP(t, \{0\})$ has $(-\log a, -\log a)$-bounded slope.

If $a=4$, we are dealing with  the Chebyshev map. This map is smoothly conjugate on the interior of $[0,1]$ to the full tent map with slope $\pm 2$, so all measures other than $\mu_0$ have Lyapunov exponent equal to $\log2$. It is straightforward to verify that there is a unique phase transition at $t = -1$. 

Now suppose $a \ne 4$. 
Any other measure in $\scrM(f_a)$ lives on $X:= [f^2(c),f(c)]$ where the norm of the derivative is bounded away from (and below) $a$. Thus $\overline{\chi}(X) < \log a$, and $\underline{\chi}(X) > -\infty$ since $a \in \scrQ$. We know $\scrP(t,X)$ has $(\overline{\chi}(X),\underline{\chi}(X))$-bounded slope.

We have $\scrP(t) = \max(\scrP(t,\{0\}), \scrP(t,X))$. Applying Lemma \ref{lem:pslope} will give a required phase transition at some parameter $t_0$. To see that $t_0 \leq -h_a /(\log a - h_a)$, 
note first that for each $a$, $h_a = \scrP(0) = \scrP(0,X)$, and that $0 = \scrP(0, \{0\})$. 
Thus if $h_a = 0$ then $t_0=0$, by uniqueness of $t_0$. 

If $h_a >0$ let $\mu_a$ be the (unique) measure of  maximal entropy. The graphs of $\scrF(t,\mu_0)$ and $\scrF(t,\mu_a)$ intersect when $t = -h_a/(4 - \chi_{\mu_a})$. The result follows upon applying Ruelle's inequality: $h_\mu \leq \max(0, \chi_\mu)$.
\eprf

\remark We shall show later that there are unimodal maps without negative Schwarzian derivative but with non-flat critical point which admit phase transitions in the negative spectrum (i.e. at some negative $t$) where the equilibrium states on both sides of the phase transition have positive entropy. 

\smallskip
\remark Pesin has conjectured that there should be a positive measure set of parameters in the neighbourhood of each Misiurewicz parameter for which, for each real $t$, there exists a unique equilibrium measure, and such that the pressure function is real analytic everywhere. The proposition above shows that one should at least restrict one's attention to measures sitting on $[f^2(c),f(c)]$ in the non-renormalisable setting. Under this additional hypothesis, we suspect the conjecture to be true. Note that near any non-renormalisable Misiurewicz parameter will be positive measure sets of both non-renormalisable and renormalisable Collet-Eckmann parameters. Indeed, such maps are accumulated by both non-renormalisable and renormalisable post-critically finite maps, around which one can apply the Benedicks-Carleson construction to find the required positive measure sets (\cite{Tsujii:PosLyap}, \cite{Thunberg:Unfolding}).

\smallskip

%

\begin{dfn} We shall call a collection of maps $\{g_\alpha\}_{\alpha \in \Delta}$ a \emph{full unimodal family} provided:
\begin{itemize}
	\item $\Delta$ is an interval;
	\item each $g_\alpha : I_\alpha \to I_\alpha$ is a $C^3$ unimodal map of the interval $I_\alpha$ and has negative Schwarzian derivative and non-flat critical point;
	\item the boundary of $I_\alpha$ depends continuously on $\alpha$;
	\item rescale $g_\alpha$ by an affine, orientation-preserving conjugacy to get a unimodal map $g^*_\alpha : [0,1] \to [0,1]$; then $\alpha \mapsto g^*_\alpha$ is continuous for the $C^1$ topology on $\{g_\alpha^*\}_{\alpha \in \Delta}$;
	\item for every $a \in (0,4]$ there is an $\alpha \in \Delta$ and a conjugacy between $g_\alpha$ and the quadratic map $f_a$.
\end{itemize}
\end{dfn}
Of course, the quadratic family $\{f_a\}_{a \in (0,4]}$ is a full unimodal family. Requiring negative Schwarzian derivative and non-flat critical point means that each periodic  attractor of $g_\alpha$ is essential (its immediate basin contains the critical point) and that wandering intervals do not exist. 

Let $\{g_\alpha\}_{\alpha \in \Delta}$ be a full unimodal family. Let $\alpha_0 \in \Delta$ be such that $g_{\alpha_0}$ is conjugate to the Chebyshev map $f_4$, but so that there are $\alpha$ arbitrarily close to $\alpha_0$ which are not conjugate to $f_4$. Then it follows from the definitions and continuity that there exists a sequence $\{\alpha_k\}_{k \geq 2}$ with $\lim_{k \to \infty} \alpha_k = \alpha_0$ such that $g_{\alpha_k}$ has critical orbit satisfying  $g_{\alpha_k}^k(c) = c$ and either 
$$
g_{\alpha_k}^i(c) < c < g_{\alpha_k}(c)
$$
or
$$
g_{\alpha_k}^i(c) > c > g_{\alpha_k}(c)
$$
for all $i$ with $2 \leq i < k$. 
The first possibility holds for maps with the same orientation as the quadratic family. Note that $c$ is contained in a restrictive interval of period $n$. For each $k \geq 2$, let $\Delta_k$ denote the connected component containing $\alpha_k$ of paramaters $\alpha \in \Delta$ such that $g_\alpha$ is renormalisable of type $k$. Denote by $J_\alpha$ the corresponding maximal restrictive intervals. If
for each $k \geq 2$, $\{g^k_\alpha :J_\alpha \to J_\alpha\}_{\alpha \in \Delta_k}$ is  a full unimodal family then we shall call $\{\alpha_k, \Delta_k\}_{k\geq 2}$  a \emph{suitable sequence}. 

The following lemma follows from 
Sections II.4-5 of
\cite{DeMeloVanStrien} (see in particular pages 148-149) which contain stronger results and details with less restrictive definitions. We only need to remark that one can choose $\alpha_k$ converging to $\alpha_0$ because the same holds for the quadratic family.
\begin{lem} \label{Fact:full} Given a full unimodal family $\{g_\alpha\}_{\alpha \in \Delta}$, there exists $\alpha_0 \in \Delta$ with $g_{\alpha_0}$  conjugate to $f_4$ and a suitable sequence $\{\alpha_k, \Delta_k\}_{k\geq2}$ with $\lim_{k\to\infty} \alpha_k = \alpha_0$. 
\end{lem}

Given any full unimodal family we fix  a suitable sequence $\{\alpha_k,\Delta_k\}_{k\geq2}$  and call maps $g_\alpha$ with $\alpha \in \Delta_k$ 
 \emph{simply renormalisable of type $k$} (these are not the only renormalisable maps of a given type $k$).  

For $\alpha \in \Delta_k$, denote by $X_\alpha$ the set of points which never enter the interior of the restrictive interval $J_\alpha$. 
Let us call a unimodal map $g:J \to J$ \emph{trivial} if $g$ only has one fixed point. In this case the orbit of every point converges to the boundary fixed point. 
\begin{lem} \label{lem:cascade}
Let $\{\alpha_k\}_{k\geq2}$ be a suitable sequence for a full unimodal family $\{g_\alpha\}_{\alpha \in \Delta}$. 
Given any $\varepsilon > 0$, there exists $K >1$ such that, for each $k\geq K$ and all $\alpha \in \Delta_k$,
$$
 1 - \varepsilon < \HD(X_\alpha).
$$
Moreover, if $\scrR g_\alpha$ is not trivial, then $\HD(X_\alpha) < 1$. 
\end{lem}
\beginpf Let $\alpha_0 := \lim_{k\to \infty} \alpha_k$. 
The map $g_{\alpha_0}$  is topologically conjugate to the Chebyshev map $f_4$. It has hyperbolic dimension equal to 1 (see e.g.\ \cite{Me:HypDim}). Therefore, denoting by $c$ the critical point of $g_{\alpha_0}$, there is a $\delta >0$ such that the set of points never entering $B(c,\delta)$ is a compact hyperbolic set of dimension greater than $1 - \varepsilon/2$. 

This set persists under perturbations, so for all $\alpha$ sufficiently close to $\alpha_0$, the set of points never entering $B(c, \delta/2)$ say, under iteration by $g_\alpha$, is a compact hyperbolic set of dimension greater than $1 - \varepsilon$. 
Taking $\alpha \in \Delta_k$ for $k$ large implies $\alpha$ is close to $\alpha_0$, by hypothesis. We now explain why the restrictive interval $J_\alpha$ is contained in $B(c,\delta/2)$ for large $k$. Suppose not, then there exists a $\nu > 0$ such that for $k$ large and $\alpha \in \Delta_k$, $J_\alpha$ compactly contains $W := B(c,\nu)$.  Some iterate of $g_{\alpha_0}$ maps $W$ onto $I_{\alpha_0}$. Thus for all $\alpha$ sufficiently close to $\alpha_0$, some iterate of $g_\alpha$ maps $W$ onto an interval compactly containing $J$, contradiction. The first statement follows. 

Now, since $g_\alpha$ has negative Schwarzian derivative, each non-repelling periodic orbit contains a critical point in its immediate basin of attraction. Thus if $\partial J_\alpha$ contains a parabolic periodic point, then $\scrR g_\alpha$ is trivial. Otherwise, $X_\alpha$ does not contain any non-repelling periodic points. Then $X_\alpha$ is hyperbolic by, for example, the Mañé hyperbolicity theorem, and so has Hausdorff dimension strictly less than one. 
\eprf

\begin{prop}\label{prop:CEpht} For each $n \geq 0$, there exists a positive measure set $A$ of parameters such that, for all $a \in A$, 
\begin{itemize}
\item 
the pressure function $\scrP_{f_a}(t)$ admits exactly $n$ phase transitions in $(0,1)$ and no additional phase transitions on neighbourhoods of 0 and of 1; 
\item
the pressure function is piecewise analytic on a neighbourhood  of $[0,1]$ and is strictly decreasing;  
\item $f_a$ admits an absolutely continuous invariant probability measure.
\end{itemize}
\end{prop}
\begin{prop} \label{prop:infrenorm} There exist uncountably many parameters $a$ for which the pressure function of the quadratic map $f_a$ admits:
\begin{itemize}
 \item (countably) infinitely many phase transitions in $(0,1)$; 
\item exactly one phase transition at some $t_a < 0$;
\item no phase transition at $t=0$  nor for $t \geq 1$.
\end{itemize}
Between the phase transitions the pressure function is analytic. For $t < 1$ it is strictly decreasing and for $t \geq 1$, $\scrP(t) = 0$. At each phase transition there exist exactly two distinct equilibrium states.
\end{prop}
\emph{Proof of Propositions \ref{prop:CEpht}, \ref{prop:infrenorm}:}
Consider sequences of integers $(a_1b_1a_2b_2\ldots a_nb_n)$ where $a_i \geq 3$ and $b_i\geq 0$. 
Denote by $\scrA_1$ the set of parameters $a \in [3,4]$ such that $f_a$ is simply renormalisable of type $a_1$. Denote by $\scrB_1 \subset \scrA_1$ the set of parameters in $\scrA_1$ such that $\scrR f_a $ is $b_1$-times renormalisable of type $2$. Define inductively $k_n$, $\scrA_n$ and $\scrB_n$ for $n>1$ as follows: 
\begin{itemize}
\item $k_1 := 0$; $k_n := k_{n-1} + 1 + b_{n-1}$;  
\item there is a subset $\Delta$ of $\scrB_n$ such that $\{\scrR^{k_n}f_a\}_{a\in\Delta}$ is a full unimodal family by Lemma \ref{Fact:full}; fix a suitable sequence and let $\scrA_{n+1}$ be the set of parameters $a \in \scrB_n$ such that $\scrR^{k_n}f_a$ is simply renormalisable of type $a_{n+1}$ and $\scrR^{k_n + 1}f_a$ is not trivial; 
\item $\scrB_{n+1}$ is the set of parameters $a \in \scrA_{n+1}$ such that $\scrR^{k_n +1}f_a$ is $b_{n+1}$-times renormalisable of type 2. 
\end{itemize}
For $n \geq 1$ set $J_n^a$ as the restrictive interval of $\scrR^{k_n}f_a$. For $n = 1$ set $L_1^a := I$ and for $n \geq 2$ define $L^a_n$ as the 
restrictive interval of $\scrR^{k_n - 1}f_a$ (so $L^a_n$ is the domain of $\scrR^{k_n}f_a$). Then define $X_n^a$ as the set of points in $L^a_n$ which never get mapped into the interior of $J_n^a$ under iteration by $\scrR^{k_{n-1}}f_a$.

Given any set $K$, we shall write $\scrO_{f_a}(K) := \bigcup_{i\geq 0}f_a^i(K)$ for the smallest, forward-invariant set for $f_a$ containing $K$.

We shall use the following inductive step.
If we fix any sequence of integers $(a_1b_1\ldots b_{n-1}a_n)$ with $a_i \geq 3$ and $b_i \geq 0$ then provided $b_n$ and $a_{n+1}$ are sufficiently large,  the following properties hold. 
\begin{enumerate}
\item
$$\ol{\chi}(\scrO_{f_a}(L_n^a)) < (1/2) \inf_{i < n} \ul{\chi}(X^a_i);$$
\item for all $i < n$
\begin{equation}  
0 < 1 - \HD(X^a_n) < \HD(X^a_n) - \HD(X^a_i).
\end{equation} 
\end{enumerate}
The first point follows on taking $b_n \geq 1$ sufficiently large. Indeed, according to the renormalisation theory of Sullivan (\cite{Sullivan:Bounds}; \cite{DeMeloVanStrien}, chapter VI), for $a \in \scrA_n$, the derivative of $\scrR^{k_n + m}f_a$ is uniformly bounded in $m$ and in $a$. Thus  the Lyapunov exponents of any $f_a$-invariant measures on $\scrO_{f_a}(L_n^a)$ are exponentially small in $b_n$.

To show the second point, we remark that, having fixed some $b_n \geq 1$,  for all $a \in \scrB_n$ there is an $\varepsilon >0$ such that $\HD(X^a_i) < 1 - \varepsilon$ for all $i < n$, by a uniform hyperbolicity argument. Then Lemma \ref{lem:cascade} gives the required $a_{n+1}$.

Thus, for each $n >0$ there are plenty of choices for a sequence $(a_1b_1\ldots a_nb_n)$ such that we have a corresponding parameter interval $\scrA_n$ and the sets $X_i^a$ verify
\begin{equation}\label{eqn:chis}
\begin{split}
\log|Df(0)| > \ol{\chi}(\scrO_{f_a}(X_i^a)) \geq \ul{\chi}(\scrO_{f_a}(X_i^a)) > (1/2) \ol{\chi}(\scrO_{f_a}(X_{i+1}^a))  \\
\end{split}
\end{equation}
and 
\begin{equation}\label{eqn:hds}
0 < 1 - \HD(X^a_{i+1}) < \HD(X^a_{i+1}) - \HD(X^a_i).
\end{equation}
for $i = 1, \ldots, n-1$ 
and all $a \in \scrA_n$.

Now apply Lemma \ref{lem:bddslope} to the restricted pressure functions $\scrP(t,\scrO_{f_a}(X^a_i))$ and $\scrP(t, \scrO_{f_a}(X^a_{i+1}))$. The slope inequalities (\ref{eqn:chis}) give a unique intersection at some $t =: t_i$. 
The dimension estimates 
(\ref{eqn:hds}), together with the slope estimates (\ref{eqn:chis}), imply by elementary geometry that 
$$
t_i < \HD(\scrO_{f_a}(X^a_i)) < t_{i+1}$$
 for each $i$, and that $\scrP(t_i) = \scrP(t_i, \scrO_{f_a}(X^a_i)) > 0$. 
Again using (\ref{eqn:chis}), and using analyticity of the restricted pressure functions, there is a phase transition at each $t_i$, $i = 1, \cdots, n$, and these are the only phase transitions on $(-\varepsilon,t_n)$ for some $\varepsilon >0$.

Let us  show Proposition \ref{prop:CEpht}. 
 Following the work of Pesin and Senti \cite{PesinSenti:Moscow}, there exists a Lebesgue positive measure subset $A$ of Collet-Eckmann parameters $a \in \scrA_n$ for which a unique equilibrium state exists for each $t$ in a neighbourhood of $[0,1]$ for $\scrR^{k_n}f_a$. 
Moreover, using the techniques of Bruin and Todd \cite{BT:EquilibriumInterval} the pressure functions for the renormalised maps for these parameters can be shown to be analytic on a neighbourhood of $[0,1]$. Thus for these parameters we also have analyticity on $(t_n, 1+\varepsilon)$ for some $\varepsilon > 0$. Since $a \in A$ are Collet-Eckmann parameters, $\ul{\chi}(I) > 0$ so the pressure function is strictly decreasing, and $f_a$ admits an absolutely continuous invariant probability measure.


For Proposition \ref{prop:infrenorm}, it suffices for us to consider infinite sequences of the form $(a_1b_1a_2b_2\ldots)$ satisfying (\ref{eqn:hds}) and (\ref{eqn:chis}) for all truncations, and the corresponding parameters $\bigcap_n \scrA_n$ (for each sequence this infinite intersection contains exactly one parameter). We remark that  $f$ restricted to $\omega(c)$, the omega-limit set of the critical point, is uniquely ergodic and its measure $\mu$ has zero entropy. The Lyapunov exponent of $\mu$ is non-negative by  \cite{Przytycki:LyapNonneg}. One can use \cite{Ledrappier:AbsCnsInterval} to show it is not strictly positive; otherwise $\mu$-almost every point would be contained in arbitrarily small (restrictive) intervals getting mapped by an iterate of $f$ onto some fixed interval, contradiction. Thus $\scrF(\mu,t) = 0$ for all $t$. Thus $\scrP(t) = 0$ for $t \geq 1$.
\eprf

\begin{prop} There exists a smooth unimodal map with non-flat critical point for which the pressure function admits phase transitions at $s$ and $t$ for some $s < 0 < t$. The pressure function $\scrP(t)$ is strictly convex and analytic on each of $(-\infty,s), (s,t), (t,+\infty)$. 
\end{prop}
\beginpf
Let $f_n$ be a quadratic map which is renormalisable of type 3, whose renormalisation is $n-1$ times renormalisable of type 2 for some $n\geq 2$, and whose final renormalisation, $\scrR^n f_n$ say, is topologically conjugate to the Chebyshev map $x \mapsto 4x(1-x)$. Let $J'$ denote the domain of $\scrR^n f_n$.

Let $X_n$ denote the largest transitive hyperbolic compact set of points which never enter the interior of the first renormalisation interval under iteration by $f_n$. Then $X_n$ does not contain $\{0\}$. Standard considerations give $H > 0$ and $\lambda > 0$ such that for all $n \geq 2$, $H < \HD(X_n)$ and $\lambda < \ul{\chi}(X_n) < \ol{\chi}(X_n) < \log 4$. 
The restricted pressure function $\scrP(t,X_n)$ is analytic, cuts the $t$-axis at $\HD(X_n) < H < 1$ and has $(-\log 4, -\lambda)$-bounded slope. We write $f$ for $f_n$ and $X$ for $X_n$ in what follows, dropping the dependence on $n$.

We want to modify $f$. 
Since $f$ is a quadratic Misiurewicz map, we have $\ul{\chi}([0,1]) >0$. 
Denote by $\beta$ the point in $\partial J'$ fixed by $\scrR^n f$, and by $\alpha$ the other, internal, fixed point of $\scrR^n f$. 
Let the open interval $V$, $\alpha \in V \subset J'$ verify the following: 
\begin{itemize}
\item $f^k(\partial V) \cap V = \emptyset$ for all $k \geq 0$;
\item $|D\phi_V| \geq 1+ \ul{\chi}([0,1])/2$, where $\phi_V$ denotes the first return map to $V$.
\end{itemize}
Note that if $h$ is some smooth function and $c_1,c_2 \in \arr$, we write $c_1 < |Dh| < c_2$ if $c_1 < |Dh(x)| < c_2$ for all $x$ in the domain of $h$.

The first point ensures that each branch of the first return map to $V$ will map its domain diffeomorphically onto $V$ (\emph{cf.}\ the \emph{nice intervals} of Martens \cite{Martens:Cantor}).
Indeed, let $A$ and $B$ be connected components of $f^{-k}(V)$ and $f^{-l}(V)$ respectively for some $k,l$ with $0 \leq k<l$ and suppose $\partial A \cap B \ne \emptyset$, so $f^l(\partial A) \cap V \ne \emptyset$. But then $f^{l-k}(\partial V) \cap V \ne \emptyset$, contradiction. This implies that the return time is constant on each connected component of the domain of the first return map. 

To find arbitrarily small $V$ satisfying the first point one can use density of periodic orbits. For the second, use negative Schwarzian derivative and  extendibility of branches onto $J'$. By the Koebe Principle (\cite{DeMeloVanStrien}, Theorem IV.1.2), taking  $V$ small enough, on each branch of $\phi_V$ the derivative is approximately constant. The lower bound on Lyapunov exponents gives a lower bound for the derivative at the fixed point of each branch (a periodic point of $f$). 

\begin{figure}[htb]
\center{\includegraphics{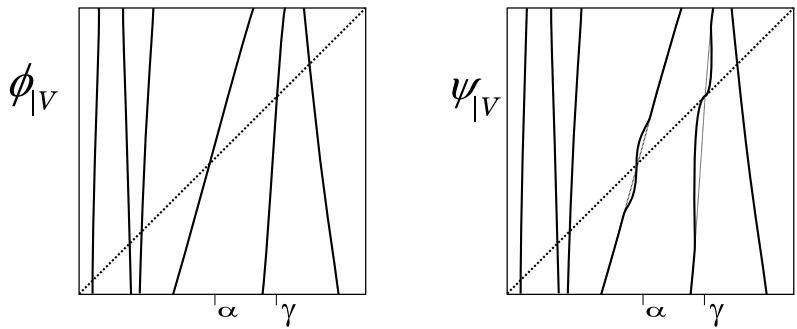}}
\caption{Modifying $f$  modifies the first return map. \label{fig:modifying}}
\end{figure}

Let $\gamma \ne \alpha$ be another fixed point of $\phi_V$ and modify $f$ (see figure \ref{fig:modifying}) on neighbourhoods of $\alpha$ and $\gamma$ compactly contained in their respective branch domains (of $\phi_V$), to get a $C^\infty$ topologically conjugate map $g$ with first return map $\psi_V$ to $V$ so that 
\begin{enumerate}
\item $|D\psi_V(\alpha)| > 2^{2^n12} $;
\item $1 < |D\psi_V(\gamma)| < 1 + \ul{\chi}_f([0,1])/2$;
\item $1+ \ul{\chi}_f([0,1])/3 \leq |D\psi_V|$.
\end{enumerate}
The first point implies that the measure sitting on the orbit of $\alpha$ (with period $2^{n-1}3$) has Lyapunov exponent greater than $4\log 4$. 
The second point implies that the measure sitting on the orbit of $\beta$ cannot be an equilibrium state (see \cite{MakStas:Nonrec} for a discussion of how this can cause a problem in the rational context).
The third point means that we have not created any parabolic or attracting orbits. 

Let $J$ denote the set $\bigcup_{i\geq 0}g^i(J')$. The restricted pressure function $\scrP_g(t,J)$ is analytic. Indeed the first return map to the interval, delimited by the two  preimages of $\alpha$ under $\scrR^n g$, is a hyperbolic induced map and the usual techniques for it (e.g. \cite{BT:EquilibriumInterval}) can be applied. All measures other than the one sitting on the orbit of $\beta$ \emph{lift} to this induced map.

\begin{figure}[htb]
\center{\includegraphics{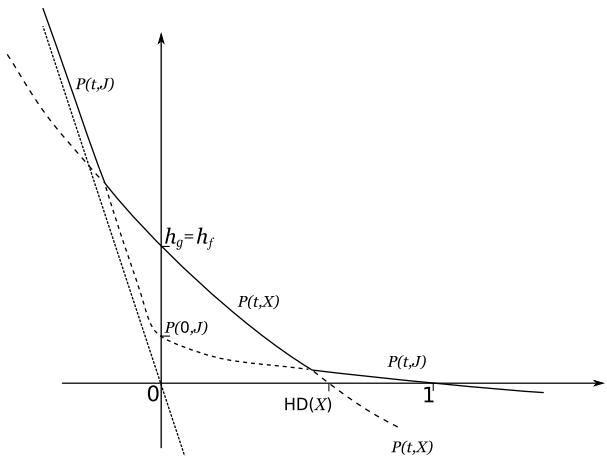}}
\caption{Graphs of the restricted pressure functions. \label{fig:posneg}}
\end{figure}

 We have $\scrP(0,J) = (\log 2)/(2^{n-1}3)$, $\scrP(1, J) = 0$. Because of convexity, for all $n$ large enough that $\scrP(0,J) < H\lambda$, this pressure function has a transverse (and unique in $(0,1)$) intersection with $\scrP(t, X)$ at some $t \in (0,1)$, see figure \ref{fig:posneg}.

Now we need to see what happens in the negative spectrum. $\scrP(t,J) \geq -t4\log4$ due to the measure sitting on the orbit of $\alpha$. The entropy $h_g = h_f$ of $g$ is greater than $(\log 2)/2$ (which is at least $2^{n-1}$ times $\scrP(0,J)$). 
One can deduce that there is an intersection of the graphs of the restricted pressure functions at some maximal $t = :s$ satisfying $-h_g/(3\log 4) < s < 0$. Then a simple calculation again gives that  the derivative
$|D\scrP(s,J)| \geq (h_g - \scrP(0,J))/|s|$. Since $n\geq 2$, 
$$
|D\scrP(s,J)| > (h_g /(2|s|) > \log 4 > |D\scrP(s,X)|$$
 so, by convexity, this is the only intersection in the negative spectrum.
\eprf

We finish with a remark which first arose in conversation with Juan Rivera-Letelier. Let $f$ be the analytic map with quadratic fixed point, renormalisable of type $k >2$, which is a fixed point of the renormalisation  operator (\emph{i.e.}\ $\scrR f$ is an affine rescaling of $f$). Then there  are exactly two phase transitions of the pressure function, one at some negative $t < 0$, the other at some $t_* > 0$ equal to the dimension of the hyperbolic set of points which never enter the interior of the restrictive interval. At $t_*$ there are an \emph{infinity} of equilibrium states, one on each level of the \emph{filtration} into transitive hyperbolic sets, and one on the omega-limit set of the critical point. 


\subsection*{Acknowledgements}  
 We would like to thank the referees and  
J.\ Rivera-Letelier and M.\ Todd for their helpful and interesting comments, questions and suggestions.